\documentclass[a4paper,fontsize=12pt,DIV=17]{scrartcl}

\usepackage[T1]{fontenc}
\usepackage[utf8]{inputenc}

\usepackage{geometry}

\usepackage{microtype}

\usepackage{fourier}
\usepackage{libertine}

\usepackage{amsmath}
\usepackage{amssymb}
\usepackage{esint}

\usepackage[utf8]{inputenc}

\newcommand{\R}{\mathbb R}
\newcommand{\N}{\mathbb N}

\newcommand{\Rset}{\mathbb{R}}
\newcommand{\Cset}{\mathbb{C}}
\newcommand{\Nset}{\mathbb{N}}
\newcommand{\Sset}{\mathbb{S}}
\newcommand{\Bset}{\mathbb{B}}
\newcommand{\Zset}{\mathbb{Z}}
\newcommand{\abs}[1]{\lvert #1 \rvert}
\newcommand{\bigabs}[1]{\bigl\lvert #1 \bigr\rvert}
\newcommand{\Bigabs}[1]{\Bigl\lvert #1 \Bigr\rvert}
\newcommand{\floor}[1]{\lfloor #1 \rfloor}
\newcommand{\dif}{\;\mathrm{d}}
\newcommand{\seminorm}[1]{\lvert #1 \rvert}
\newcommand{\norm}[1]{\lVert #1 \rVert}

\newcommand{\st}{\ \big\vert \ }

\usepackage[pdftitle={Sobolev mappings: from liquid crystals to irrigation via degree theory},
pdfauthor={Jean Van Schaftingen}]{hyperref}

\begin{document}

\author{Jean \textsc{Van Schaftingen} \thanks{Universit\'{e} catholique de Louvain, Institut de Recherche en Mathématique et Physique, Chemin du Cyclotron 2 bte L7.01.01, 1348 Louvain-la-Neuve, Belgium.}}
\title{Sobolev mappings: from liquid crystals to irrigation via degree theory}
\date{\small February 3, 2017} 

\maketitle

\thispagestyle{empty}
 
\begin{abstract}
Sobolev spaces are a natural framework for the analysis of problems in partial differential equations and calculus of variations. Some physical and geometric contexts, such as liquid crystals models and harmonic maps, lead to consider Sobolev maps, that is, Sobolev vector functions 
whose range is constrained in a surface or submanifold of the space.
This additional nonlinear constraint provokes the appearance of finite-energy topological
singularities. These singularities are characterized by a nontrivial topological invariant such as the topological degree, they represent an obstruction to the strong approximation by smooth maps and they become source and sink terms in an optimal transportation or irrigation problem of topological charges arising in the study of the weak approximation and of the relaxed energy.
\end{abstract}

\vspace{\stretch{1}}
These notes are based on the \emph{Godeaux lecture} that I gave during the \emph{Brussels Summer School in Mathematics} in Brussels on Wednesday August~3, 2016. 
I thank the Brussels Summer School of Mathematics and the Belgian Mathematical Society for their invitation to give this lecture.

I have kept this text close to the original oral and blackboard presentation. 
In particular the reference list is limited to an extent that is reasonable for an oral presentation. I invite the interested reader to start from the bibliographies of these references to dive into the abundant relevant publications by a whole community of mathematicians.

\newpage

\section{Classical Sobolev spaces}

Before considering Sobolev spaces between manifolds, we revisit the classical theory of Sobolev spaces of scalar functions.

\subsection{Motivation}

One motivation to the study of Sobolev spaces is \emph{Laplace's equation}
\begin{equation}
\label{equationLaplace}
 \left\{
   \begin{aligned}
     -\Delta u & = 0 && \text{in \(\Omega\)},\\
     u & = g && \text{on \(\partial \Omega\)}.
   \end{aligned}
 \right.
\end{equation}
Here the set \(\Omega\) is a \emph{bounded open set} of the Euclidean space \(\Rset^m\) whose \emph{boundary \(\partial \Omega\) is a smooth submanifold} of \(\Rset^m\), or equivalently, \(\partial \Omega\) is a compact set of \(\Rset^m\) which is locally and up to an affine isometry of \(\Rset^m\) the graph of smooth function on a subset of \(\Rset^{m - 1}\); a typical example would be an \(m\)--dimensional ball.
The operator \(\Delta\) is the \emph{Laplacian} (Laplace's operator), which can be defined for each \(x \in \Omega\) by
\[
 \Delta u (x) = \operatorname{div} \nabla u (x)
 = \operatorname{tr} D^2 u (x)
 = 2 (m + 2) \lim_{r \to 0} \frac{\fint_{B_r (x)} u  - u(x)}{r^2},
\]
that is,
\[
 \fint_{B_r (x)} u
  = u (x) + \frac{r^2}{2 (m + 2)}\, \Delta u (x) + o (r^2),
\]
as \(r \to 0\).
The Laplacian can thus be interpreted geometrically as the coefficient of the \emph{asymptotic deviation of the average value} of the function on a ball from its value at the center of the ball. 
In the one-dimensional case \(m = 1\), the Laplacian evaluates the \emph{convexity} of the function; in higher-dimensions it is related to the \emph{mean curvature} of the graph. 

Laplace's equation \eqref{equationLaplace} arises in Fourier's linear model of the stationary temperature distribution in a homogeneous medium with prescribed temperature at the boundary (\emph{steady-state heat equation});
it also appears in steady state \emph{diffusion} models (Fick's law), 
in \emph{electrostatics} (electric potential in the absence of charges) and in computation of expected values of a \emph{Brownian motion} when it hits the boundary. 

If we assume that we have two functions \(u \in C^2 (\Omega) \)  and \(v \in C^1 (\Omega)\), then we can compute by expansion of the scalar product
\[
 \int_{\Omega} \abs{D v}^2
 = \int_{\Omega} \abs{D u}^2
 + 2 \int_{\Omega} D u \cdot D (v - u)
 + \int_{\Omega} \abs{D (v - u)}^2.
\]
By integration by parts (\emph{Gauß--Ostrogradsky divergence theorem}), we have 
\[
 \int_{\partial \Omega} (v - u)\, \partial_{\nu} u
 = \int_{\Omega} \operatorname{div}\, \bigl((v - u) \nabla u \bigr)
 = \int_{\Omega} D (v - u) \cdot D u
 + \int_{\Omega} (v - u)\, \Delta u,
\]
(here \(\nu\) denotes the exterior normal vector to the submanifold \(\partial \Omega\))
and therefore, if the function \(u\) is a solution of Laplace's equation \eqref{equationLaplace} in the domain \(\Omega\) and if \(u = v\) everywhere on the boundary \(\partial \Omega\), then 
\[
 \int_{\Omega} D u \cdot D  (v - u) = 0.
\]
In particular, we have 
\[
 \int_{\Omega} \abs{D v}^2
 = \int_{\Omega} \abs{D u}^2
 + \int_{\Omega} \abs{D (v - u)}^2
 \ge \int_{\Omega} \abs{D u}^2,
\]
with equality if and only if \(D v = D u\) in \(\Omega\).
This suggests searching for a solution \(u\) by minimizing the 
\emph{Dirichlet energy functional}
\[
  \int_{\Omega} \abs{D v}^2
\]
among the functions \(v \in C^1 (\Omega) \cap C (\Bar{\Omega})\) that satisfy the boundary condition \(v = g\) on \(\partial \Omega\). 

If we want to actually \emph{construct a solution}, we first consider a sequence of functions \((u_k)_{k \in \Nset}\)
in \(C^1 (\Omega) \cap C (\Bar{\Omega})\) such that 
\[
  \lim_{k \to \infty} \int_{\Omega} \abs{D u_k}^2 =
  c =
  \inf \,\Biggl\{\int_{\Omega} \abs{D v}^2 \st v \in C^1 (\Omega) \cap C (\Bar{\Omega}) \Biggr\}
  \ge 0.
\]
We have, for each \(k, \ell \in \Nset\), the parallelogram identity 
\[
 \int_{\Omega} \abs{D u_k - D u_\ell}^2
 + \int_{\Omega} \abs{D u_k + D u_\ell}^2
 = 2 \int_{\Omega} \abs{D u_k}^2 + 2 \int_{\Omega} \abs{D u_\ell}^2.
\]
We also observe that 
\[
 \int_{\Omega} \abs{D u_k + D u_\ell}^2
 = 4 \int_{\Omega} \abs{D w_{k, \ell}}^2 \ge 4 c,
\]
where \(w_{k, \ell} = \frac{u_k + u_\ell}{2}\), and \(w_{k, \ell} = g\) on \(\partial \Omega\).
Therefore 
\[
 \lim_{k, \ell \to \infty}\int_{\Omega} \abs{D u_k - D u_\ell}^2 = 0. 
\]
This is a \emph{Cauchy condition} for the sequence \((u_k)_{k \in \Nset}\) in a certain
\emph{seminorm}:
\[
  \seminorm{u}_{W^{1, 2}} 
  = \norm{D u}_{L^2} 
  = \Bigl(\int_{\Omega} \abs{D u}^2 \Bigr)^\frac{1}{2}.
\]

There are two issues. First, the quantity \(\seminorm{\cdot}_{W^{1, 2}}\) is \emph{not positive definite} on functions: it vanishes on constant functions. Next, the set \(C^1 (\Omega) \cap C (\Bar{\Omega})\) is \emph{not complete} under this norm, even up to a quotient by constant functions.

To tackle the first problem, we consider instead the norm 
\[
 \norm{u}_{W^{1, 2}} = \Bigl(\int_{\Omega} \abs{u}^2 + \abs{D u}^2 \Bigr)^\frac{1}{2}.
\]
It can be proved by a clever integration by parts
that, under the assumption that the set \(\Omega\) is bounded in one direction or has finite measure, there exists a constant \(C \in \Rset\) such that for
every 
\(u \in C^1 (\Omega) \cap C (\Bar{\Omega})\) such that \(u = 0\) on \(\partial \Omega\), 
we have \emph{Poincaré's inequality}
\[
 \int_{\Omega} \abs{u}^2 \le C \int_{\Omega} \abs{D u}^2.
\]
In the language of spectral theory, the spectrum of the Laplacian \(\Delta\) with Dirichlet boundary conditions on such a domain \(\Omega\) has a negative spectrum that stays away from \(0\).
In particular, for each \(k, \ell \in \Nset\) we have, since \(u_k - u_\ell = 0\) on \(\partial \Omega\),
\[
 \int_{\Omega} \abs{u_k - u_\ell}^2
 + \int_{\Omega} \abs{D u_k - D u_\ell}^2
 \le (C + 1)\int_{\Omega} \abs{D u_k - D u_\ell}^2,
\] 
and thus 
\[
 \lim_{k, \ell \to \infty} \int_{\Omega} \abs{u_k - u_\ell}^2
 + \int_{\Omega} \abs{D u_k - D u_\ell}^2 = 0.
\]
Hence, our sequence is a Cauchy sequence for this apparently stronger semi-norm.

To deal with the second issue, we consider the Sobolev space \(W^{1, 2}(\Omega)\) defined as the \emph{completion of the space} 
\[
  \bigl\{\, u\in C^1 (\Omega) \st \norm{u}_{W^{1, 2}} < +\infty\,\bigr\},
\]
under the norm \(\norm{\cdot}_{W^{1, 2}}\).
That is, we consider the smallest normed vector space \(W^{1, 2} (\Omega)\) that contains the above space and such that any Cauchy sequence for this Sobolev norm converges with respect to the same Sobolev norm. This space can be obtained by quotienting the space of Cauchy sequences by its linear subspace of sequences converging to \(0\).
In particular on this space the Dirichlet functional will achieve a minimum that will be a
\emph{weak solution} of Laplace's equation \eqref{equationLaplace}; 
\emph{regularity theory} can show in many cases that this solution is in fact a classical solution.

The Sobolev  space \(W^{1, 2} (\Omega)\) is part of the \emph{larger family of Sobolev spaces} obtained by considering the norm
\[
 \norm{u}_{W^{1, p}} = \Biggl(\int_{\Omega} \abs{u}^p + \abs{D u}^p \Biggr)^\frac{1}{p},
\]
and the Sobolev space \(W^{1, p} (\Omega)\) defined as the completion of the set
\[
  \bigl\{ \,u \in C^1 (\Omega) \st \norm{u}_{W^{1, p}} < +\infty \,\bigr\},
\]
under the norm \(\norm{\cdot}_{W^{1, p}}\).
This generalization will allow us to \emph{observe high-dimensional phenomena in low dimensions}. 
It is also relevant for example in regularity theory in the calculus of variations and partial differential equations because it provides intermediate levels of regularity of solutions.

\subsection{Properties of Sobolev spaces}

The Sobolev spaces can be characterized in different ways
\begin{itemize}
 \item functions that are absolutely continuous on almost every line,
 \item functions in \(L^p\) that have a distributional gradient that can be represented by an \(L^p\) vector field,
 \item when \(\Omega\) is a cube and \(p = 2\), they have a Fourier characterization:
 \[
  \norm{u}_{W^{1, 2}}^2 = \sum_{k \in \Zset^m} (1 + \abs{k}^2)\, \abs{\widehat{u} (k)}^2,
 \]
 where the function \(\widehat{u} : \Zset^m \to \Cset\) gives the multidimensional Fourier coefficients \(\widehat{u} (k)\) of \(u\) for \(k \in \Zset^m\).
\end{itemize}

The Sobolev spaces are naturally embedded into the space \(L^p (\Omega)\), and thus they can be identified with equivalence classes of functions. 

In classical analysis, it is well-known that differentiability of a function at a point implies continuity at that point .
This is still the case for Sobolev spaces in the high-integrability case \(p > m\) where every Sobolev function can be identified with a continuous function. Indeed, the \emph{Morrey inequality} states that there exists a constant \(C > 0\) such that if \(x, y \in \Omega\) and if \(u \in C^1 (\Omega) \cap W^{1, p}( \Omega)\), then 
\begin{equation}
\label{equationMorrey}
 \abs{u (x) - u (y)} \le C\, \Biggl(\int_{\Omega} \abs{D u}^p \Biggr)^\frac{1}{p} \abs{x - y}^{1 - \frac{m}{p}}.
\end{equation}
In particular, any function \(u \in W^{1, p} (\Omega)\) has a continuous representative.
This representative turns out to be differentiable almost everywhere in \(\Omega\).
The exponent in the Morrey inequality \eqref{equationMorrey} is optimal: if we define the function \(u : \Bset^m\to \R\) for \(x \in \Bset^m \setminus \{0\}\) by
\[
 u (x) = \abs{x}^{\gamma},
\]
with \(\gamma > 1 - \frac{m}{p}\), then \(u \in W^{1, p} (\Bset^m)\).

When \(p < m\), the Sobolev embedding theorem states that there exists a constant \(C > 0\) such that for every function \(u \in C^1 (\Omega) \cap W^{1, p}( \Omega)\),
\[
 \Bigl(\int_{\Omega} \abs{u}^{\frac{m p}{m - p}}\Bigr)^{1 - \frac{p}{m}}
 \le C \int_{\Omega} \abs{D u}^p.
\]
The result is optimal: if we define the function \(u : \Bset^m\to \R\) for each \(x \in \Bset^m \setminus \{0\}\) by
\[
 u (x) = \frac{1}{\abs{x}^{\gamma}},
\]
with \(\gamma > \frac{m}{p} - 1\), then \(u \in W^{1, p} (\Bset^m)\).
Although they are not continuous, Sobolev functions      are defined on the boundary 
by the trace inequality which states that for every function \(u \in C (\Bar{\Omega}) \cap W^{1, p} (\Omega)\),
\[
 \Bigl(\int_{\partial \Omega} \abs{u}^\frac{(m - 1)p}{m - p} \Bigr)^{1 - \frac{p - 1}{m - 1}} \le C \int_{\Omega} \abs{D u}^p .
\]
More generally Sobolev functions are well-defined on every \(k\)--dimensional submanifold when \(k > m - p\). Moreover, the case \(p = 1\) and \(k = m -1\) is an additional exceptional case where the trace is well-defined.

In the \emph{critical case} \(p = m > 1\), Sobolev functions are \emph{neither continuous nor bounded}. Indeed, if the function \(u : \Bset^m \to \R\) is defined for each \(x \in \Bset^m \setminus \{0\}\) by 
\[
 u (x) = \Bigl(\log \frac{1}{\abs{x}} \Bigr)^\gamma,
\]
with \(0 < \gamma < 1 - \frac{1}{m}\), then \(u \in W^{1, m} (\Bset^m)\).
Functions in \(W^{1, m} (\Rset^m)\) are in all the \(L^p (\Rset^m)\) spaces and even better: they are exponentially integrable (Moser--Trudinger inequality) and have vanishing mean oscillation (VMO).

\section{Sobolev mappings}

\subsection{Liquid crystals and definition}

In a simplified model, \emph{liquid crystal molecules} have an orientation at every point that can be represented by a \emph{unit vector}, that is, an element of the two-dimensional sphere \(\Sset^2\) \cite{Mermin1979,BallZarnescu2011,Brezis1991}. 
(For many liquid crystal, in fact only the direction matters and not the orientation and liquid crystal are then modelled by a distribution of lines, or equivalently, a map into the real projective plane \(\Rset P^2\).)

Similarly to the temperature, that tends to minimize variations, the simplest model states that the \emph{liquid crystal minimizes the Dirichlet energy functional }
\[
 \int_{\Omega} \abs{D u}^2,
\]
where \(\abs{Du}\) is the Hilbert--Schmidt or Frobenius norm given by 
\[
  \abs{D u}^2 = \sum_{i, j = 1}^m \abs{\partial_i u^j}^2.
\]
Such a minimizing map is a \emph{harmonic map}, which generalizes the notion of geodesic: when \(m = 1\), harmonic maps are in fact constant velocity parametrizations of geodesics.

A \emph{natural space} to consider such problems is the \emph{nonlinear Sobolev space} of maps from \(\Omega\) to the manifold \(\Sset^{n - 1}\):
\[
 W^{1, p} (\Omega, \Sset^{n - 1})
 = \,\bigl\{ \, u \in W^{1, p} (\Omega, \Rset^n) \st u \in \Sset^{n - 1} \text{ almost everywhere in \(\Omega\)} \,\bigr\},
\]
defined by constraining vector Sobolev functions to take their value into the sphere \(\Sset^{n - 1}\).

\subsection{Density problem}

An alternative way to define the Sobolev space \(W^{1, p} (\Omega, \Sset^{n - 1})\) would have been to take the closure 
of the set \(C^1 (\Omega, \Sset^{n - 1})\) with respect to the \(W^{1, p}\)--norm which induces naturally a distance on this set.
This leads to the question whether the definitions coincide.

\paragraph{Supercritical case}
When \(p > m = \dim \Omega\), if we consider a sequence \((u_k)_{k \in \Nset}\) in \(C^1 (\Omega, \Rset^n)\) that approximates the Sobolev mapping \(u \in W^{1, p} (\Omega, \Sset^{n - 1})\)
in the space \(W^{1, p} (\Omega, \Rset^{n})\), we have as a consequence of the Morrey inequality \eqref{equationMorrey}, for almost every \(x \in \Omega\),
\[
\begin{split}
 \abs{u_k (x) - u (x)} &\le 
   \fint_{\Omega} \abs{u_k (x) - u (x) - (u_k (y) - u (y))} \dif y
   + \fint_{\Omega} \abs{u_k (y) - u (y)} \dif y\\
 &\le C\, \norm{u_k - u}_{W^{1, p}}.
\end{split}
\]
For almost every \(x \in \Omega\), we have in particular, since \(\abs{u (x)} = 1\),
\[
 \abs{u_k (x) - 1} \le \abs{u_k (x) - u (x)} \le C \norm{u_k - u}_{W^{1, p}} \to 0.
\]
This implies that for \(k \in \Nset\) 
large enough, for almost every \(x \in \Omega\), \(\abs{u_k (x)} \ge \frac{1}{2}\), and thus the function
\[
 \Tilde{u}_k= \frac{u_k}{\abs{u_k}}
\]
is well-defined on the set \(\Omega\) and takes its value into \(\Sset^{n - 1}\).
Moreover, we have 
\[
 D \Tilde{u}_k = P_{u_k}\frac{D u_k}{\abs{u_k}},
\]
where \(P_{u_k}\) denotes the projection orthogonal to \(u_k\), which is well defined almost everywhere on \(\Omega\).
Up to a subsequence \(\abs{D u_k} \le g\) and \(u_k \to u\) almost everywhere in \(\Omega\), for some \(g \in L^p (\Omega)\).
In particular we have, almost everywhere in \(\Omega\) 
\[
  \abs{D \Tilde{u}_k} \le \frac{\abs{D u_k}}{\abs{u_k}} \le 2 g, 
\]
and, as \(k \to \infty\),
\[
 D \Tilde{u}_k \to D u,
\]
almost everywhere in \(\Omega\)
and therefore, by Lebesgue's dominated convergence theorem,
\[
 \int_{\Omega} \abs{D \Tilde{u}_k - D u}^p \to 0.
\]
Similary, we have \(\Tilde{u}_k \to u\) almost everywhere in \(\Omega\) as \(k \to \infty\),
and for each \(k \in \Nset\), \(\abs{\Tilde{u}_k} = 1\)
almost everywhere on \(\Omega\),
so that 
\[
 \int_{\Omega} \abs{\Tilde{u}_k - u}^p \to 0.
\]
Essentially, smooth maps are dense in \(W^{1, p} (\Omega, \Sset^n)\) because every Sobolev map in \(W^{1, p}\) is continuous.

\paragraph{Critical case}
The approach fails if \(p = m = \dim \Omega\), because there is no Morrey inequality anymore in this case.
In fact a map \(u \in W^{1, m} (\Omega, \Sset^{n - 1})\) need not be continuous: consider for example the map \(u : \Bset^m \setminus \{0\} \to \Sset^1\) defined for \(m \ge 2\) and each \(x \in \Bset^{m} \setminus \{0\}\) by
\[
 u (x) = \biggl( \cos \Bigl(\log \tfrac{1}{\abs{x}}\Bigr)^\gamma, \sin \Bigl(\log \tfrac{1}{\abs{x}}\Bigr)^\gamma\biggr),
\]
with \(\gamma < 1 - \frac{1}{m}\). Since the map \(u\) is not continuous it cannot be approximated uniformly by smooth maps, or even by continuous ones, because the uniform convergence preserves continuity.

Although these Sobolev maps are not continuous, smooth maps are still dense \cite{SchoenUhlenbeck1983}.
Sobolev maps can be approximated if we use the standard approximation by averaging with convolution: 
since the domain \(\Omega\) is smooth, the function can be extended by a reflection with respect to the boundary to a function \(\Tilde{u} \in W^{1, m} (\Tilde{\Omega}, \Sset^{n - 1})\) such that \(\Tilde{\Omega}\) is open and \(\Tilde{\Omega} \supset \Bar{\Omega}\).
We then approximate the extended map \(\Tilde{u}\) by averages: when \(\delta > 0\) is small enough so that 
\(B_{\delta} (x) \subset \Omega\), we take
\[
 u_{\delta}
 (x)
 =\fint_{B_\delta (x)} \Tilde{u}.
\]
It can be proved that
\[
 \lim_{\delta \to 0} \;\norm{u_{\delta} - u}_{W^{1, m}} \to 0
\]
Moreover, for almost every \(x \in \Omega\),
\[
\begin{split}
 \bigabs{\abs{u_\delta(x)}  - 1}
 &= \fint_{B_\delta (x)} \abs{u_\delta (x) - \Tilde{u} (y)}\dif y
 = \fint_{B_\delta (x)} \Bigabs{\fint_{B_\delta (x)} \Tilde{u} - \Tilde{u} (y)}\dif y\\
 &\le  \fint_{B_\delta (x)} \fint_{B_\delta (x)} \abs{ \Tilde{u} (z) - \Tilde{u} (y)}\dif z \dif y
 \le C\, \Bigl(\int_{B_\delta (x)} \abs{D \Tilde{u}}^m\Bigr)^\frac{1}{m},
 \end{split}
\]
where the last step is a variant of the classical \emph{Poincaré inequality}.
Therefore, we have 
\[
 \lim_{\delta \to 0} \operatorname*{ess\, sup}_{x \in \Omega} \,\bigabs{\abs{u_\delta (x)} - 1} = 0.
\]
We can thus continue as in the supercritical case.

The essential ingredient in this proof is that critical Sobolev maps have vanishing mean oscillation and that the latter property is sufficient for the strong approximation \cite{BrezisNirenberg1995,BrezisNirenberg1996}.

\paragraph{The almost critical case: topological degree}
When \(n \le p < n + 1 \le m\), smooth maps are not dense in \(W^{1, p} (\Omega, \Sset^n)\)  \cite{BethuelZheng1988}.
Let us look at this when \(m = n = 2\) and \(\Omega = \Bset^2\) is the unit ball in \(\Rset^2\).
We consider the map \(u : \Bset^2 \to \Sset^1\) defined for each \(x \in \Bset^2 \setminus \{0\}\) by
\[
 u (x) = \frac{x}{\abs{x}}.
\]
We compute for every \(x \in \setminus \{0\}\)
\begin{align*}
  D u (x) &= \frac{\abs{x}^2 \operatorname{id} - x \otimes x}{\abs{x}^3}&
  &\text{and thus }&
  \abs{D u (x)} &= \frac{1}{\abs{x}}.
\end{align*}
One has 
\[ 
  \int_{\Bset^2} \abs{Du}^p = \int_{\Bset^2} \frac{1}{\abs{x}^p} \dif x = 2 \pi \int_0^1 \frac{1}{r^{p - 1}} \dif r < +\infty
\]
if and only if \(p < 2\).
We observe, that by the change of variable to polar coordinates and Fubini's theorem, we have if the function \(f : \Bset^2 \to \R\) is integrable
\[
 \int_{\Bset^2} f
 = \int_0^1 \int_{\Sset^1} f (r \cos \theta, r \sin \theta) \dif \theta r \dif r. 
\]
Therefore, if \((u_k)_{k \in \Nset}\) is an approximating sequence of the map \(u\),
we have, up to a subsequence, for almost every \(r \in (0, 1)\),
\[
 \lim_{r \to 0} \int_0^{2 \pi} \abs{u_{k, r}' - u_r'}^p + \abs{u_{k,r} - u_r}^ p\to 0, 
\]
where the function \(u_{k, r}: [0, 2\pi] \to \Rset\) is defined for each \(\theta \in [0, 2 \pi]\) by
\[
 u_{k, r} (\theta) = u_k (r \cos \theta, r \sin \theta).
\]
(We are dropping the radial part of the derivative.)

We now make the following observation: for every \(r \in (0, 1)\) and \(k \in \Nset\), we have 
\[
 \int_0^{2 \pi} u_{k, r} (\theta) \wedge u_{k, r}' (\theta)\dif \theta = 0.
\]
Indeed, by differentiation and integration by parts,  
\[
\begin{split}
 \frac{d}{dr} \int_0^{2 \pi} u_{k, r} (\theta) \wedge u_{k, r}' (\theta)\dif \theta 
 &= \int_0^{2 \pi} \partial_r u_{k, r} (\theta) \wedge u_{k, r}' (\theta) \dif \theta\\
 &\qquad+ \int_0^{2 \pi}  u_{k, r} (\theta) \wedge \partial_r u_{k, r}' (\theta) \dif \theta\\
 &= 2 \int_{0}^{2 \pi}  \partial_r u_{k, r} (\theta) \wedge u_{k, r}' (\theta) \dif \theta = 0,
\end{split}
\]
because both vectors \(\partial_r u_{k, r} (\theta)\) and \(u_{k, r}' (\theta)\) are tangent to the unit circle at \(u_{k, r}(\theta)\) and thus colinear.
Moreover,
\[
 \lim_{r \to 0} \int_0^{2 \pi} u_{k, r} (\theta) \wedge u_{k, r}' (\theta)\dif \theta 
 = 0.
\]
Therefore, for every \(r \in (0, 1)\) and \(k \in \N\)
\[
 \int_0^{2 \pi} u_{k, r} (\theta) \wedge u_{k, r}' (\theta)\dif \theta = 0.
\]
By the approximation assumption, we have
\[
  \int_0^{2 \pi} u_{r} (\theta) \wedge u_{r}' (\theta) \dif \theta = 
  \lim_{k \to \infty}  \int_0^{2 \pi} u_{k, r} (\theta) \wedge u_{k, r}' (\theta)\dif \theta = 0
\]
where for every \(r \in (0, 1)\) and \(\theta \in [0, 2 \pi]\),
\[
 u_r (\theta) = u (r \cos \theta, r \sin \theta) = (\cos \theta, \sin \theta),
\]
so that 
\[
 \int_0^{2 \pi} u_{r} (\theta) \wedge u_{r}' (\theta) \dif \theta= 2 \pi.
\]

The quantity 
\[
 \frac{1}{2\pi} \int_0^{2 \pi} u_{r} (\theta) \wedge u_{r}' (\theta) \dif \theta,
\]
defines in fact the \emph{winding number} or \emph{degree} \(\deg u_r\) of the map \(u_r : \Sset^1 \to \Sset^1\), which takes only integral values: if we write \(u_r (\theta)
= (\cos \phi (\theta), \sin \phi (\theta))\), then 
\[
 \int_0^{2 \pi} u_{r} (\theta) \wedge u_{r}' (\theta) \dif \theta
 = \int_0^{2 \pi} \phi' (\theta) \dif \theta = \phi (2 \pi) - \phi (0) \in 2 \pi \Zset.
\]
The degree also provides a way of computing the \emph{winding number} for a \emph{continuous map}.
In fact, it allows to classify the connected components of the set \(C (\Sset^1, \Sset^1)\) of continuous maps from the circle to the circle.
The formula of the \emph{Kronecker index}
\[
 \fint_{\Sset^{m - 1}} \det (D w)
\]
defines a similar object for a map \(w \in C^1 (\Sset^{m - 1}, \Sset^{m - 1})\).
Brouwer and Hadamard have showed how this allows generalised this contruction to the classification of connected components of \(C (\Sset^{m - 1}, \Sset^{m - 1})\) (\emph{Brouwer topological degree}).
This allows one to show the same obsruction for \(W^{1, p} (\Bset^m, \Sset^{m - 1})\). 

In this case maps that are smooth except at a finite number of points are dense in \(W^{1, p} (\Omega, \Sset^{m - 1})\) \cite{BethuelZheng1988}.
The density of maps that are smooth except at a finite number of points can be interpreted by saying that Sobolev maps are topologically like continuous maps defined on a smaller set obtained by drilling holes in the original domain.

More generally we can consider the space of Sobolev maps \(W^{1, p} (\Omega, N)\), for which there is an obstruction to the density of smooth maps if the homotopy group \(\pi_{\floor{p}} (N)\) is nontrivial, that is, if there exists a continuous map on the \(\floor{p}\)--dimensional sphere \(f \in C(\Sset^{\floor{p}}, N)\) which is not the restriction of a continuous map on the \(\floor{p + 1}\)--dimensional ball \(\Bar{f} \in C(\Bset^{\floor{p +1 }}, N)\) \cite{BethuelZheng1988}. Here \(\floor{p}\) denotes the integer part of the real number \(p\).
This is the only restriction on the target space when the domain \(\Omega\) has a trivial topology \cite{Bethuel1991,HangLin2003}.
In general, maps from \(\Omega\) to \(N\) that are smooth outside a \((m - \floor{p})\)--dimensional set are dense in the Sobolev space \(W^{1, p}(\Omega, N)\).

\section{Connecting singularities and optimal transport}

For a map \(u \in W^{1, m}(\Omega, \Sset^{n - 1})\), we define the \emph{relaxed energy} to be 
\begin{multline*}
 \mathcal{E} (u) = \inf\; \Bigl\{ \liminf_{k \to \infty} \int_{\Omega} \abs{D u_k}^m \st u_k \to u \text{ almost everywhere in \(\Omega\)}\\
 \text{ and \(u_k \in C(\Omega, \Sset^{n - 1}) \cap W^{1, m}(\Omega, \Sset^{n - 1})\)}\Bigl\}.
\end{multline*}
The relaxed energy \(\mathcal{E}\) can be interpreted as an alternate definition of the energy by density; it is in fact the lower semicontinuous relaxation of the energy defined on smooth functions.
It is related to the energy by a Fatou-like property
\[
 \mathcal{E} (u) \ge \int_{\Omega} \abs{D u}^m,
\]
with equality 
if and only if \(u\) can be approximated in \(W^{1, m}\) by smooth functions.

Let us compute this quantity for maps which are smooth outside a finitely set of singularity points; such maps are dense in the space \(W^{1, m} (\Omega, \Sset^{n -1}) \).
This computation was originally performed in the case \(m = n = 3\), motivated by liquid crystal models 
\cite{BrezisCoronLieb1986,BethuelBrezisCoron1990};
we will consider the geometrically simple case \(m = n = 2\) and \(p = 1\) \cite{BrezisMironescuPonce2005}.

First, we consider the case of a zero-degree homogeneous function on the ball \(B_R\) for which \(u (x) = u (x/\abs{x})\), with a singularity at the origin \(0 \in \R^2\). In the \emph{topologically trivial} case \(\deg u \vert_{\partial B_R} = 0\), the map \(u\) can be lifted as \(u = (\cos \varphi, \sin \varphi)\) with a zero-degree homogeneous function \(\varphi : B_R \to \R\).
We then take a smooth extension of the restriction of the lifting \(\varphi \vert_{\partial B_R}\) to \(B_R\) and a homogeneous extension in \(\Rset^2\setminus B_R\), resulting in a smooth function \(\Tilde{\varphi} : \R^2 \to \R\). We then set for \(\delta > 0\) and each \(x \in B_R\)
\[
 u_\delta (x) = \bigl(\cos \varphi (x/\delta), \sin \varphi (x/\delta)\bigr).
\]
and observe that 
\[
  \int_{B_R} \abs{D u_\delta} = \int_{B_1 \setminus B_\delta} \abs{D u}
  + \frac{C}{\delta} \int_{B_\delta} \abs{D \Tilde{\varphi}(x/\delta)}  \dif x
  \le \int_{B_R} \abs{D u} + \frac{\delta}{R} \int_{B_R} \abs{D \varphi},
\]
which implies by letting \(\delta \to 0\) that the relaxed energy coincides with the energy
\[
  \mathcal{E} (u)= \int_{B_R} \abs{D u}.
\]
(The reader will observe that in fact the family \((u_\delta)_{\delta > 0})\) converges strongly to \(u\) in \(W^{1, 1} (B_R, \Sset^1)\) as \(\delta \to 0\).)

In the case \(\deg u\vert_{\partial B_R} = k \ne 0\) where \(u\) has a \emph{nontrivial topological singularity}, we write 
\(u (r \cos \theta, r \sin \theta) = (\cos \varphi (\theta), \sin \varphi (\theta))\),
with a function \(\varphi \in C^1 ([0, 2 \pi])\) such that \(\varphi (2 \pi) = \varphi (0) + 2 k \pi\).
We now define the function \(u_\eta : B_R \to \Sset^1\) by setting for each \(\theta \in [0, 2 \pi]\),
\begin{multline*}
 u_\eta (r \cos \theta, r \sin \theta)\\
 =
 \left\{
 \begin{aligned}
   &\bigl(\cos \varphi (\tfrac{\theta}{1 - \eta}), \sin \varphi (\tfrac{\theta}{1 - \eta})\bigr) & &\text{if \(0 \le \theta \le (1 - \eta) 2 \pi\),}\\[1em]
   &\bigl(\cos \bigl(\tfrac{1 - \theta/(2 \pi)}{\eta} \varphi (2 \pi) + 
   \tfrac{\theta/(2 \pi) + \eta - 1}{\eta} \varphi (0)\bigr),\\
   &\qquad
   \sin \bigl(\tfrac{1 - \theta/(2 \pi)}{\eta} \varphi (2 \pi) + 
   \tfrac{\theta/(2 \pi) + \eta - 1}{\eta} \varphi (0)\bigr)
   \bigr) && \text{if \((1 - \eta) 2 \pi \le \theta \le 2 \pi\).}
 \end{aligned}
 \right.
\end{multline*}
The two key observations are that 
\[
 \int_{B_R} \abs{D u_\eta} = \int_{B_R} \abs{D u} + \abs{k}\, 2 \pi R.
\]
and that \(u_\eta\) has degree \(0\) and is homogeneous, so that we can conclude that 
\[
 \mathcal{E} (u) = \int_{B_R} \abs{D u} + \abs{k}\, 2 \pi R.
\]
This also works if the map \(u\) is only homogeneous in a neighbourhood of the singularity: the additional radial derivative terms will become negligible as \(\eta \to 0\).

If we now have \emph{two topological singularities} \(a_1\) and \(a_2\) in \(\Rset^2\) around which the map \(u\) is homogeneous and has respectively degree \(1\) and \(-1\), we have, by decomposing the domain in suitable tangent balls and performing the constructions on these balls
\[
  \mathcal{E} (u) = \int_{\Rset^2} \abs{D u} + 2 \pi\, \abs{a_1 - a_2}.
\]
In the case with \emph{four topological singularities} \(a_1\), \(a_2\), \(a_3\) and \(a_4\), of degrees \(1\), \(1\), \(-1\) and \(-1\), we have 
\[
 \mathcal{E} (u) = \int_{\Bset^2} \abs{D u} + 2 \pi\, \min\,\bigl\{\,\abs{a_1 - a_3} + \abs{a_2 - a_4}, \;
 \abs{a_1 - a_4} + \abs{a_2 - a_3}\,\bigr\}.
\]
In general, we have singularities \(a_1, \dotsc, a_N\) of degree \(d_1, \dotsc, d_N\). If \(\int_{\Rset^2} \abs{Du} < +\infty\), then  we have automatically \(d_1 + \dotsb + d_N = 0\).
We then have 
\[
  \mathcal{E} (u) = \int_{\Bset^2} \abs{D u} + 2 \pi \,\inf\;\biggl\{\, \sum_{i, j = 1}^N b_{ij}\, \abs{a_i - a_j}
 \; \st\; \sum_{i = 1}^N b_{ij} - b_{ji} = d_j \text{ and } b_{ij} \in \Nset \,\biggr\}.
\]
It can be observed that
\begin{multline*}
\inf\;\biggl\{\, \sum_{i, j = 1}^N b_{ij}\, \abs{a_i - a_j}
 \; \st\; \sum_{i = 1}^N b_{ij} - b_{ji} = d_j \text{ and } b_{ij} \in \Nset \,\biggr\}\\
= \inf\;\biggl\{\, \sum_{i, j = 1}^N b_{ij}\, \abs{a_i - a_j}
 \; \st\; \sum_{i = 1}^N b_{ij} - b_{ji} = d_j \text{ and } b_{ij} \ge 0 \,\biggr\},
\end{multline*}
because the second infimum always has an integral optimum \cite[\S 3.3]{Wolsey1998}.

This infimum can be thought as the solution of an \emph{optimal transport problem}. The topological singularities with positive degree can be thought as offering some degree and the topological singularities with negative degree as demanding some degree, and the problem is to match the offer and demand at the least transportation cost. 
The \emph{Kantorovich duality} is an important tool for such transport problems.
The idea is that instead of having a transporter charging for transport, 
we could have an economically equivalent situation of a trader buying all the degree where it is offered and reselling it where it is demanded. Instead of having producers and consumers trying to minimize the transport cost, we would have the trader trying to maximize its profit. 
But, the trader cannot have his business working if the producers and consumers can manage to transport at a cost lower than then margins. It turns then out that the minimal transport cost is always going to be the maximal trade profit.

In our particular case, the duality for the transport problem states that 
\begin{multline*}
 \inf\;\biggl\{ \,\sum_{i, j = 1}^N b_{ij}\, \abs{a_i - a_j}
  \st \sum_{i = 1}^N b_{ij} - b_{ji} = d_j \text{ and } b_{ij} \ge 0 \,\biggr\}\\[-1em]
  = \sup\; \biggl\{ \, \sum_{i = 1}^N f_i d_i \st f_i - f_j \le \abs{a_i - a_j}\, \biggr\},
\end{multline*}
where \(f_i\) is the price of a singularity at \(a_i\).
In the language of functions, we can define a function \(\varphi: \{a_1, \dotsc, a_N\} \to \Rset\) such that \(\varphi (a_i) = f_i\), and we have thus 
\begin{multline*}
 \inf\;\biggl\{ \,\sum_{i, j = 1}^N b_{ij}\, \abs{a_i - a_j} \st \sum_{i = 1}^N b_{ij} - b_{ji} = d_j \text{ and } b_{ij} \ge 0 \,\biggr\}\\[-1em]
  = \sup\; \biggl\{\, \sum_{i = 1}^N \varphi (a_i)\, d_i \st \varphi : \{a_1, \dotsc, a_N\} \to \Rset
   \text{ and }\varphi (a_i) - \varphi (a_j) \le \abs{a_i - a_j} \,\biggr\}.
\end{multline*}
The condition on the function \(\varphi\) is merely that it is \(1\)--Lipschitz on the finite set \(\{a_1, \dotsc, a_N\}\).  Any such function is a restriction of a \(1\)--Lipschitz function on the entire space \(\Rset^2\). This follows from Kirszbraun's extension theorem or can be seen directly here by setting for each \(x \in \Rset^2\)
\[
  \varphi (x) = \min\,\bigl\{\,\varphi (a_i) + \abs{x - a_i} \st i \in \{1, \dotsc, N\}\, \bigr\}. 
\]
We can thus rewrite the quantity as 
\begin{multline*}
 \inf\;\biggl\{\, \sum_{i, j = 1}^N b_{ij}\, \abs{a_i - a_j} \st \sum_{i = 1}^N b_{ij} - b_{ji} = d_j \text{ and } b_{ij} \ge 0 \,\biggr\}\\[-1em]
 =\sup\; \biggl\{ \,\sum_{i = 1}^N \varphi (a_i)\, d_i \st \varphi :\Rset^2 \to \Rset
  \text{ and }\abs{\varphi (x) - \varphi (y)} \le \abs{x - y} \,\biggr\}.
\end{multline*}
Interestingly, the right-hand side of this last identity can be computed without knowing the place of the singularities \(a_1, \dotsc, a_N\).
Indeed one has
\[
\begin{split}
 \sum_{i = 1}^m \varphi (a_i)\, d_i
 &= \int_{\Rset^2} (u \wedge  Du) \wedge D \varphi\\
 &= \int_{\Rset^2} u_1 \,\partial_1\! u_2 \,\partial_2\! \varphi - u_2 \,\partial_1\! u_1 \,\partial_2\! \varphi
 - u_1 \,\partial_2\! u_2 \,\partial_1\! \varphi + u_2 \,\partial_2\! u_1 \,\partial_1\! \varphi.
\end{split}
\]
This makes sense if \(\varphi\) is continuously differentiable with a bounded derivative but also when \(\varphi\) is merely Lipschitz. We have thus 
\begin{multline*}
 \inf\;\Bigl\{ \sum_{i, j = 1}^N b_{ij}\, \abs{a_i - a_j}
  \st \sum_{i = 1}^N b_{ij} - b_{ji} = d_j \text{ and } b_{ij} \ge 0 \Bigr\}\\
  = \sup\; \Bigl\{ \int_{\Rset^2} (u \wedge  Du) \wedge D \varphi \st \varphi :\Rset^2 \to \Rset \text{ and }\abs{\varphi (x) - \varphi (y)} \le \abs{x - y} \,\Bigr\}.
\end{multline*}
The right-hand side makes sense for every weakly differentiable map \(u :\Rset^2 \to \Sset^1\) such that \(\int_{\Rset^2} \abs{D u} < + \infty\) and can be used to compute the relaxed energy \(\mathcal{E} (u)\) for any such map.

\section{Sublinear energies and irrigation}

The topological singularities  for maps in \(W^{1, 3} (\Rset^4, \Sset^2)\), 
are related to the homotopy classes of \(C (\Sset^3, \Sset^2)\) which are classified by the \emph{Hopf degree}. 
By performing similar computations on the relaxed energy, for a map \(u\) having
two singularities of Hopf degree \(d\) at the points \(a_1\) and \(a_2\), we have 
\[
 \mathcal{E} (u) = \int_{\Omega} \abs{D u}^3 + c (d)\, \abs{a_1 - a_2},
\]
with 
\[
 c (d) \simeq \abs{d}^{\frac{3}{4}}.
\]
The striking difference with the previous case is the appearance of the \emph{sublinear exponent} \(\frac{3}{4}\), in connection with the \emph{Hopf fibration} \cite{HardtRiviere2003}.

This \emph{sublinear cost of transport} will change the picture of the optimal transport of singularities as soon as one has four singularities \(a_1\), \(a_2\) of four singularities \(a_1\), \(a_2\), \(a_3\) and \(a_4\), of degrees \(1\), \(1\), \(-1\) and \(-1\), we have to consider the quantity 
\begin{multline*}
 \min\,\bigl\{\,c (1)\, \abs{a_1 - a_3} + c (1)\, \abs{a_2 - a_4},
 c (1)\, \abs{a_1 - a_4} + c (1)\, \abs{a_2 - a_3}, \\
 \,\min_{b_1, b_2 \in \R^4} c (1)\, \abs{a_1 - b_1} + c (1)\, \abs{a_2 - b_1} 
 + c (2)\, \abs{b_2 - b_1}\\[-1em]
 + c (1)\, \abs{b_2 - a_3} + c (1)\, \abs{b_2 - a_4}
 \,\bigr\}.
\end{multline*}
That is, we also need to consider the possibilities of the two singularities of degree \(1\) connecting into a degree \(2\) singularity that propagates and then splits to connect to the two singularities of degree \(-1\). 
This happens if \(c (2) < 2 c (1)\) and if \(\abs{a_1 - a_2}+ \abs{a_3 - a_4}\) is significantly smaller than \(\abs{a_1 - a_3}\). 
The sublinearity of the transport costs makes it interesting to merge singularities together at some \emph{branching point}. 

This phenomenon is related to economies of scale in networks, when increasing the cost of segment is sublinear with respect to its capacity, as in public transport networks in which people accept to take a longer path or distribution of water and electricity in which the consumers do not have individual connection up to the supplier.

We will now examine what happens when many singularities are placed in an array and a cost of transport \(c (d) = \abs{d}^\alpha\).
To understand the dimensionality issues, we will work in a Euclidean space \(\Rset^m\) and we consider
an array of \(2^{m n}\) of points source \(d\) which are at mutual distance and we consider the cost of connecting them to a sink of intensity \(-2^{mn}d\) at the centre of the cube .

A first way to do this would be to \emph{connect all the sources directly to the sink, resulting in a cost} 
\begin{equation}
\label{centralizedPlan}
 \sum_{k \in \{-2^n + \frac{1}{2}, -2^{n - 1} + \frac{3}{2}, \dotsc, 2^{n - 1} + \frac{1}{2}\}^m}
 \abs{k h} \, d^\alpha
 \simeq \frac{d^\alpha}{h^m} \int\limits_{[-2^{n - 1} h , 2^{n - 1} h]^m}  \abs{x} \dif x
 \simeq d^\alpha 2^{n (m+1)} h .
\end{equation}

Another strategy is to \emph{compute hierarchically} a cost \(C(n)\) by induction over \(n\).
When \(n = 0\) the sink and the source coincide and thus 
\[
 C (0) = 0.
\]
If \(n \ge 1\), then we can consider \(2^m\) subarrays of \(2^{m (n - 1)}\) sources and connecting them to a collecting point at the center of the array. By induction hypothesis the cost of connecting to the sources to the connecting point is \(C (n - 1)\); since they are \(2^m\), the total cost of connection between sources and collecting points in this arrangement is \(2^m C (n - 1)\). Then we connect all the collecting points to the central sink. The length of each connection is then going to be \(2^{n - 2} h \sqrt m\) and each of them has a capacity \(2^{m (n - 1)}\).
Therefore we have
\[
\begin{split}
 C (n) & = 2^{m} \bigl(C (n - 1) +  2^{n - 2} h \sqrt{m}\, (2^{m(n - 1)}d)^{\alpha} \bigr)\\
 & = 2^{m} C (n - 1) + 2^{m -1}  \sqrt{m}\, d^\alpha h\,  2^{ (n - 1) (m \alpha +1)}.
\end{split}
\]
The solution of this recurrence equation is given for \(n \in \Nset\) by 
\[
 C (n)
 = 
 \left\{
\begin{aligned}
  &\frac{2^{m - 1} \sqrt{m}\, d^\alpha h}{2^{m \alpha + 1} - 2^m}
  \,\bigl(2^{(m \alpha + 1)n} - 2^{m n}\bigr) & & \text{if \(\alpha \ne 1 - \tfrac{1}{m}\)},\\
  &  \frac{\sqrt{m}\, d^\alpha h}{2}
  2^{m n} n & & \text{if \(\alpha = 1 - \tfrac{1}{m}\)}.
\end{aligned}
 \right.
\]
The dominating term in this formula changes when \(\alpha\) crosses the \emph{critical value} \(1 - \frac{1}{m}\), and we have thus 
\[
 C (n)
 \simeq  
 \left\{
  \begin{aligned}
  & d^\alpha h \,2^{m n} & & \text{if \(\alpha < 1 - \tfrac{1}{m}\)},\\
  & d^{1 - \frac{1}{m}} h \,2^{m n} n & & \text{if \(\alpha = 1 - \tfrac{1}{m}\)},\\
  & d^\alpha h \,2^{(m \alpha + 1)n} & & \text{if \(\alpha > 1 - \tfrac{1}{m}\)}.
  \end{aligned}
\right.
\]
For \(\alpha \in (0, 1)\), this hierarchical transport plan gives an asymptotically better upper bound than the one given by the centralized plan in \eqref{centralizedPlan}; when \(\alpha =1\) both approaches are asymptotically equivalent.
It turns out that \(C (n)\) is asymptoticaly the optimal costs \cite{Bethuel}.

A particularly interesting régime is when \(h = \frac{1}{2^n}\) and \(d = \frac{1}{2^{mn}}\), where as \(n \to \infty\), the sources converge to the \(m\)--dimensional Lebesgue measure on a unit cube. One has then:
\[
  C (n)
 \simeq  
 \left\{
  \begin{aligned}
  & \frac{1}{h^{m (1 - \alpha) - 1}} & & \text{if \(\alpha < 1 - \tfrac{1}{m}\)},\\
  & \ln h & & \text{if \(\alpha = 1 - \tfrac{1}{m}\)},\\
  & 1 & & \text{if \(\alpha > 1 - \tfrac{1}{m}\)}.
  \end{aligned}
\right.
\]
The upper bound is bounded if and only if \(\alpha > 1 - \frac{1}{m}\), and thus \emph{the Lebesgue measure is irrigable if and only if \(\alpha > 1 - \frac{1}{m}\)} \cite{DevillanovaSolimini2007}.

Going back to the problem of the relaxed energy in \(W^{1, 3} (\Rset^4, \Sset^2)\),
it appears that since \(m = 4\) the exponent \(\alpha = \frac{3}{4} = 1 - \frac{1}{m}\) is critical. This means that connecting \(2^{4 n}\) singularities of Hopf degree \(1\) on an array will carry a cost \(2^{3 n} n\). On the other hand it is possible \cite{Bethuel} to construct a map \(u_n \in W^{1, 3} (\Rset^4, \Sset^3)\) having such a set of singularities with 
\[
 \int_{\Rset^4} \abs{D u_n}^3 = O (2^{3 n}).
\]
We have thus 
\[
 \lim_{n \to \infty} \frac{\mathcal{E} (u_n)}{\displaystyle \int_{\Rset^4} \abs{D u_n}^3} = +\infty.
\]

In linear functional analysis, the classical Banach--Steinhaus \emph{uniform boundedness principle} states that if \((L_n)_{n \in \Nset}\) is a sequence of bounded linear operators on a Banach space \(X\) and if \((x_n)_{n \in \Nset}\) is a sequence in \(X\) such that \(\norm{L_n x_n}/\norm{x_n} \to +\infty\) as \(n \to \infty\), then there exists \(x \in X\) such that \(\norm{L_n x} \to +\infty\) as \(n \to \infty\). This can be restated by defining the quantity \(\mathcal{E} (x) = \sup_{n \in \Nset} \norm{L_n x} \in [0, +\infty]\) and saying that if the function \(\mathcal{E} (x)/\norm{x}\) is not bounded on \(X\), then there exists \(x \in X\) such that \(\mathcal{E} (x)= +\infty\).

Sobolev spaces of mappings are not linear spaces and the approximation scheme behind the relaxed energy has no reason to be linear, however there is still a \emph{nonlinear uniform boundedness} for the weak-bounded approximation problem \cite{HangLin2003III} that there exists a Sobolev mapping \(u \in W^{1, 3}(\Bset^4, \Sset^2)\) such that \(\mathcal{E} (u) = +\infty\). This particular function \(u\) shows how an infinite collection whose total energy is finite can require an \emph{infinite energy to be approximated}.



\end{document}